\documentclass[11pt]{article}
\usepackage{amsmath,amsthm}

\def\refeq#1{\if\workingver y(\ref{#1})-[[#1]]\else(\ref{#1})\fi}
\def\refth#1{\if\workingver y\ref{#1}-[[#1]]\else\ref{#1}\fi}
\def\mylabel#1{\if\workingver y\label{#1}{\bf\ \ [[#1]]\ \ }
\else\label{#1}\fi}
\def\mybibitem#1{\if\workingver y\bibitem{#1}{\bf\ \ [[#1]]\ \ }
\else\bibitem{#1}\fi}

\newtheorem{thm}{Theorem}
\newtheorem{lem}[thm]{Lemma}

\newtheorem{cor}[thm]{Corollary}

\newtheorem{algo}[thm]{Algorithm}
\newtheorem{rem}[thm]{Remark}
\newtheorem{conj}[thm]{Conjecture}

\def\begeq#1{\begin{equation}\mylabel{#1}}
\def\endeq{\end{equation}}

\def\begalg{\begin{alg}}
\def\endalg{\end{alg}}

\newcounter{li}

\def\begalg#1{\begin{algo}\mylabel{#1}\normalshape:\\}

\let\workingver=n

\def\al{{\alpha}}

\begin{document}

\title{\bf  $p^q$-Catalan Numbers and Squarefree Binomial Coefficients}
\author{Pantelimon St\u anic\u a
\thanks{On leave from the Institute of
Mathematics of Romanian Academy, Bucharest,  Romania}\\
\small Auburn University Montgomery, Department of Mathematics\\
\small Montgomery, AL 36117, USA,\\
\small\em (334) 244-3321, e-mail:  stanpan@strudel.aum.edu
}
\date{\today}
\maketitle

\vspace{2cm}

{\bf Running Head:}

$p^q$-Catalan Numbers

\newpage
\begin{abstract}
\baselineskip=2\baselineskip
In this paper we consider the generalized Catalan numbers
 $F(s,n)=\frac{1}{(s-1)n+1}\binom{sn}{n}$, which we call
$s$-Catalan numbers. We find all natural numbers $n$ such that
for $p$\, prime,
$p^q$ divides $F(p^q,n),\, q\geq 1$ and all
distinct residues of $F(p^q,n) \pmod {p^q}$, $q=1,2$. As a byproduct
we settle a
question of Hough and the late Simion on the divisibility of
the $4$-Catalan numbers by $4$.
We also prove that $\binom{p^qn+1}{n}$, $p^q\leq 99999$, is squarefree
for $n$ sufficiently large (explicit), and with
the help of the generalized Catalan numbers we find the set of
possible exceptions. As consequences, we obtain that $\binom{4n+1}{n}$,
$\binom{9n+1}{n}$ are squarefree for $n\geq 2^{1518}$, respectively
$n\geq 3^{956}$, with at most $2^{18.2}$, respectively $3^{15.3}$
possible exceptions.
\end{abstract}

\vspace{3cm}


\newpage

\title{\bf \Large $p^q$-Catalan Numbers and Squarefree Binomial Coefficients}
\maketitle

\pagestyle{myheadings}

\baselineskip=1.82\baselineskip

\section{Introduction}

Problems involving binomial coefficients were considered by
many mathematicians for over two centuries.
 R.K. Guy in \cite{G} mentions several
problems on divisibility of binomial coefficients (see {\bf B31,
B33}).  Erd\"os conjectured that for
$n>4$, $\binom{2n}{n}$ is never squarefree. This was proved by
S\'ark\"ozy in \cite{Sar}, for sufficiently large $n$, and by
Granville and Ramar\'e in \cite{GR} for any $n>4$.

Many people (see, for instance,
\cite{BS, Chu, HP,  AMM, K, Sands, Stanley}) proposed and studied
the following generalization of classical Catalan numbers
$\frac{1}{n+1}\binom{2n}{n}$,
 which we will call
{\em $s$-Catalan numbers}, $F(s,n)=\frac{1}{(s-1)n+1}\binom{s
n}{n}$. There are many interpretations of this sequence
\cite{Chu, HP, K, Sands, Stanley}, for instance:
the
number of $s$-ary trees with $n$ source-nodes, the number of ways
of associating $n$ applications of a given $s$-ary operator, the
number of ways of dividing a convex polygon into $n$ disjoint
$(s+1)$-gons with nonintersecting diagonals, and the number of
$s$-good paths (below the line $y=sx$) from $(0,-1)$ to $(n,(s-1)n-1)$.

Naturally,  some of the questions proposed by Erd\"os on
the classical Catalan numbers, may be
asked here as well, as Hough and the late Simion  proposed \cite{AMM}:
(a) When $p$ is prime, for what values of $n$ is $F(p,n)$
divisible by $p$? (b)$^*$ For what values of $n$ is $F(4,n)$
divisible by 4? (c)$^*$ What can you say when $s$ takes on the
other composite values?
There are no known answers for $(b),(c)$.
In this paper we give a simple proof to (a),  and we show that
$F(p^2,n)$ is divisible by $p^2$, unless $(p^2-1)n+1$ is an even power
of $p$, or a sum of odd powers of $p$ (with
the numbers of distinct powers summing to $p$), thereby proving (b),
and (c) for $s=p^2$.
We also prove that $\binom{p^qn+1}{n}$, $p^q\leq 99999$, is squarefree
for $n$ sufficiently large (explicit), and with
the help of the generalized Catalan numbers we find the set of
possible exceptions. As consequences, we obtain that $\binom{4n+1}{n}$,
$\binom{9n+1}{n}$ are squarefree for $n\geq 2^{1518}$, respectively
$n\geq 3^{956}$, with at most $2^{18.2}$, respectively $3^{15.3}$
possible exceptions.

\section{Preliminary Results}

Let $[x]$ be the largest integer smaller than $x$.
In this section we state a few results which will be needed later.
Lucas (1878) (see \cite{D}) found
a simple method to find $\binom{m}{n}\pmod p$.
\begin{thm}[Lucas]
If $p$ is prime, then
\(
\binom{m}{n}\equiv \binom{[m/p]}{[n/p]}\binom{m_0}{n_0}\pmod p,
\)
where $m_0,\, n_0$ are the least non-negative residues modulo $p$
 of $m$, respectively $n$.
\end{thm}

Define $n!_p$ to be the product of all integers $\leq n$,
that are not divisible by $p$. We see that
$\displaystyle n!_p=\frac{n!}{[n/p]! p^{[n/p]}}$.
Granville in \cite{AG} proves the following beautiful
generalization of Lucas' Theorem.
\begin{thm}[Granville]
Suppose that the prime power $p^q$ and positive integers\/
$m=n+r$ are given.
Let $N_j$ be the least positive residue of $[n/p^j]\pmod {p^q}$
for each $j\geq 0$
(that is, $N_j=n_j+n_{j+1} p+\cdots +n_{j+q-1} p^{q-1}$): also
make the corresponding
definitions for $m_j,M_j,r_j,R_j$. Let $e_j$ be the number of
indices $i\geq j$ for
which $m_i< n_i$ (that is, the number of {\em carries}, when
adding $n$ and $r$
in base $p$, on or beyond the $j$th digit). Then
\[
\frac{1}{p^{e_0}} \binom{m}{n}\equiv (\pm 1)^{e_{q-1}}
\frac{M_0!_p}{N_0!_p\, R_0!_p}
\frac{M_1!_p}{N_1!_p\, R_1!_p}
\cdots \frac{M_d!_p}{N_d!_p\, R_d!_p}\pmod {p^q},
\]
where $(\pm)$ is $(-1)$ except if $p=2$ and $q\geq 3$.
\end{thm}

In 1808 Legendre showed that the exact power of $p$ dividing $n\,!$ is
\begin{equation}
\label{legendre}
[n/p]+[n/p^2]+[n/p^3]+\cdots.
\end{equation}
We define (see \cite{AG}) the {\em sum of digits function}
$\sigma_p(n)=n_0+n_1+\cdots +n_d$, if $n=n_0+n_1 p+\cdots + n_d p^d$.
Then, using $\sigma$, \refeq{legendre} transforms into
\begin{equation}
\label{pow_fact}
\frac{n-\sigma_p(n)}{p-1}.
\end{equation}
We will need the following result which belongs to Kummer
\begin{thm}[Kummer]
The power to which the prime $p$ divides the binomial
coefficient $\binom{m}{n}$,
say ${v_p}(\binom{m}{n})$,
is given by the number of  carries when we add $n$ and
$m-n$ in base $p$.
\end{thm}

Our first result gives a complete answer to the first posed question (a),
generalizing the well-known result on Catalan numbers, or equivalently,
on middle binomial  coefficients
(see \cite{G}), which states
that $4\,|\, \binom{2n}{n}$, unless $n=2^k$.
\begin{thm}
Let $p$ be a prime. Then,
$p$ divides $F(p,n)$,
 unless $n$ is of the form $\frac{p^k-1}{p-1},\ k\in {\bf N}$, in which
 case $F(p,n)\equiv 1\pmod p$.
\end{thm}

\begin{proof}
We re-write $\displaystyle F(p,n)=\frac{1}{(p-1)n+1}\binom{pn}{n}=
\frac{1}{pn+1} \binom{pn+1}{n}$.
We shall find the values
$n$ such that $p\not| F(p,n)$. Since $p\not| F(p,0)$, we assume $n\not=0$.
 Applying Lucas' Theorem
repeatedly for the base $p$
representations ($0\leq m_i,n_i\leq p-1$),
$m=m_0+m_1 p+\cdots +m_d p^d$ and $n=n_0+n_1 p+\cdots +n_d p^d$, we obtain
\(
\binom{m}{n}\equiv \binom{m_0}{n_0}\binom{m_1}{n_1}\cdots
\binom{m_d}{n_d}\pmod p.
\)
For $m=pn+1\equiv 1\pmod p$, we get
\[
F(p,n)\equiv \binom{pn+1}{n}\equiv \binom{1}{n_0}\binom{n_0}{n_1}\cdots
\binom{n_{d-1}}{n_d}\pmod p,\quad n_d\not=0.
\]
Let $n_{-1}=1$.
Since for any $j$, $n_j<p$, if $\binom{n_{i-1}}{n_i}\not= 0$, then
$p\not|\binom{n_{i-1}}{n_i}$. Thus, the numbers $n$ with
$p\,\not|\, F(p,n)$ are those
positive integers $n$ with $\binom{n_{i-1}}{n_i}\not= 0,
i=0,1,\ldots,d,$ or $n=0$.
From $n_{i-1}\geq n_i$, $n_{-1}=1,n_d\not=0$, and
$\binom{n_{i-1}}{n_i}\not= 0$,
we obtain
$n_{i-1}=n_i$. Thus, $n=1+p+\cdots +p^d=\frac{p^{d+1}-1}{p-1}$.
Therefore, if $p$ is prime, $p|F(p,n)$ for any number
$n\not=\frac{p^k-1}{p-1},k\in {\bf N}$. Using Lucas' Theorem we easily see
that for $n=\frac{p^k-1}{p-1}$, the least residue of
$F(p,n)$ modulo $p$ is $1$.
\end{proof}
The following lemma will be extensively used throughout the paper
\begin{lem}
\label{pow_lem}
We have
\begin{eqnarray*}
\label{pow_pq0}
{v_p}(F(p^q,n))=
\frac{\sigma_p((p^q-1)n+1)-1}{p-1}.
\end{eqnarray*}
\end{lem}
\begin{proof}
We use the identity $F(p^q,n)=\frac{1}{p^qn+1}\binom{p^qn+1}{n}$.
Using \refeq{pow_fact} we get that the power of $p$ dividing
$\binom{m}{n}$ is
\begin{equation}
\label{pow}
{v_p}\left(\binom{m}{n}\right)=\frac{\sigma_p(n)+
\sigma_p(m-n)-\sigma_p(m)}{p-1}.
\end{equation}
Let $m=p^qn+1$. Thus, \refeq{pow} becomes
\[
\displaystyle {v_p}(F(p^q,n))=\frac{\sigma_p(n)+\sigma_p((p^q-1)n+1)-
\sigma_p(p^q n+1)}{p-1}=
\frac{\sigma_p((p^q-1)n+1)-1}{p-1},
\]
since $\sigma_p(p^q n+1)=\sigma_p(n)+1$.
\end{proof}

\section{Scarce squarefree $p^2$-Catalan numbers}

Denote by $n=(ab\ldots)_p$ the base $p$ representation of $n$,
$a$ being the most significant bit.
Our next result refers to the third question of
Hough and Simion.
\begin{thm}
\label{thm7}
Given a prime $p$, $p^2$ divides $F(p^2,n)$,
unless\/ $n$ is of the form
$\displaystyle\frac{p^{2t}-1}{p^2-1},\ t\in {\bf N}$, in which case
$F(p^2,n)\equiv 1\pmod {p^2}$, or of the form
$\displaystyle\frac{c_1 p^{2i_1+1}+\cdots +c_s p^{2i_s+1}-1}{p^2-1},\
i_1< \cdots< i_s,$
with $\displaystyle\sum_{i=1}^s c_i=p$, $s\in {\bf N}$, $0\leq c_i<p$,
in which case
$F(p^2,n)\equiv  \displaystyle \binom{p}{c_1,c_2,\ldots,c_s} \pmod{p^2}$
(the multinomial coefficient).
\end{thm}

\begin{proof}
As before $F(p^2,n)=\frac{1}{p^2n+1}\binom{p^2n+1}{n}$.
A number $n$, which does not satisfy the divisibility, must satisfy
(see Lemma \refth{pow_lem})
\begin{eqnarray}
\label{pow_p}
{v_p}(F(p^2,n))=\frac{\sigma_p((p^2-1)n+1)-1}{p-1}\leq 1,
\end{eqnarray}
which implies $\sigma_p((p^2-1)n+1)\leq p$.

Assume first that $\sigma_p((p^2-1)n+1)=1$.  Therefore,
$(p^2-1)n+1=p^k\equiv (-1)^k\pmod {p^2-1}$,
therefore $k$ must be even, say $k=2t$,
so
\(\displaystyle
n=\frac{p^{2t}-1}{p^2-1}.
\)

\def\al{{\alpha}}
Assume now that $\sigma_p((p^2-1)n+1)=l$, and $1<l\leq p$. It follows
that $(p^2-1)n+1=p^{\al_1}+\cdots +p^{\al_l}$,
$\al_1\leq \al_2\leq \cdots \leq\al_l$.
Therefore,
\(
(p^2-1)n+1\equiv l\equiv 1\pmod {p-1},
\)
and since $1<l \leq p$, we get that $l=p$. Then,
$(p^2-1)n+1=p^{\al_1}+\cdots + p^{\al_p}$,
$\al_1\leq \al_2\leq\cdots \leq \al_p$.
It follows that
\begin{eqnarray*}
(p^2-1)n
& \equiv & -1+ p\sum_{\al_i\ \text{odd}} 1 +  \sum_{\al_i\ \text{even}}1\\
& \equiv & -1+p^2-(p-1)\sum_{\al_i\ \text{even}} 1 \\
& \equiv & -(p-1)\sum_{\al_i\
\text{even}}1 \pmod {p^2-1},
\end{eqnarray*}
so $\displaystyle \sum_{\al_i\ \text{even}}1$ must be divisible by $p+1$.
Since $\displaystyle 0\leq \sum_{\al_i\ \text{even}}1\leq p$,
we see that
$\displaystyle \sum_{\al_i\ \text{even}} 1$ must be
an empty sum. Therefore, all $\al_j=2i_j+1$.
We obtain $\displaystyle n=\frac{p^{2i_1+1}+\cdots+p^{2i_p+1}-1}{p^2-1},$
$i_1\leq i_2\leq \cdots\leq i_p$, and the first claim is proved.

 Let $n_{-1}=0$.
Consider $\displaystyle n=\frac{p^{2t}-1}{p^2-1}$. It follows that
$n=(1010\cdots 101)_p$\ and $p^2 n+1$ attaches to this string the
block $01$ to the right, so it is of the same form.
Since $M_i=N_{i-2}$ and $R_i=1$, except for $R_{2t-1}=p$, we get,
using Granville's theorem,
\begin{equation}
\label{resid_p}
\frac{1}{(p^2-1)n+1}\binom{p^2 n}{n}\equiv
\frac{1}{p^2 n+1}\binom{p^2 n+1}{n}\equiv p^{e_0}
(-1)^{e_1} \frac{M_0!_p M_1!_p}{R_{2t-1}!_p}\pmod {p^2}.
\end{equation}
Now, $M_0=m_0+m_1 p=m_0+ n_{-1}p=1$,
$M_1=m_1+m_2 p=n_{-1}+n_0 p=p$ and $R_{2t-1}=p$,
implies $M_0!_p=1, M_1!_p=p!_p=(p-1)!$ and
$R_{2t-1}!_p=(p-1)!$. Thus, \refeq{resid_p}  becomes
$F(p^2,n)\equiv p^{e_0} (-1)^{e_1}\equiv 1\pmod {p^2}$, since $e_0=e_1=0$.

Consider $\displaystyle n=\frac{ c_1 p^{2i_1+1}+\cdots +
c_s p^{2i_s+1}-1}{p^2-1}$,  $i_1< \cdots< i_s$ and $c_1+c_2+\cdots +c_s=p$.
Observe that $s\geq 2$.
It follows that
\begin{eqnarray}
n&=& \frac{ c_s(p^{2i_s+1}-p)+\cdots +c_1(p^{2i_1+1}-p)+p^2-1}{p^2-1}\\
 &=& c_s(p^{2i_s-1}+p^{2i_s-3}+\cdots+1)+\cdots +
 c_1(p^{2i_1-1}+p^{2i_1-3}+\cdots+1)+1\\
 &=& c_s p^{2i_s-1}+c_s p^{2i_s-3}+\cdots +
 (c_s+c_{s-1}) p^{2 i_{s-1}-1} +\cdots +1.
\end{eqnarray}
But $p^2 n+1$ attaches the block $01$ to the right of the base $p$
representation of $n$, and since there is a carry in this case,
we get $e_0=e_1=1$.
Also, $n_0=1$, $M_0!_p=1, M_1!_p=(p-1)!$,
$R_i!_p=1$ except for $\displaystyle R_{2i_k}!_p=
(c_k p)!_p=\frac{(c_kp)!}{c_k! p^{c_k}}$ and
$R_{2i_k+1}!_p= (c_k)!_p=c_k!$, for $k=1,2,\ldots,s$.
Now, applying Granville's theorem we get
\begin{equation}
\label{congr_pq}
\begin{split}
F(p^2,n)
& \equiv  p^{e_0} (-1)^{e_1} \frac{M_0!_p M_1!_p}{ R_0!_p\cdots
R_{d+2}!_p} \equiv  (-1) p\ \frac{(p-1)!}{\displaystyle
\prod_{k} (c_k p)!_p\, c_k!}\\
&\equiv  (-1) p\ \frac{(p-1)!}{
\prod_{k} \frac{(c_k p)!}{p^{c_k}}}
 \equiv
\frac{(-1) p^{p+1} (p-1)!}{\displaystyle\strut\prod_k (c_k p)!}\pmod {p^2},
\end{split}
\end{equation}
since $\displaystyle \sum_{k=1}^s c_k=p$.
We prove that the last expression is the multinomial coefficient
(this was observed by one of our referees, whom we thank).
First, assume $p=2$. Since $s=2$ in this case, we get $c_1=c_2=1$ and the
claim is trivially satisfied. Let $p>2$. We observe that
\[
\frac{(mp+1)\cdots (mp+p-1)}{(p-1)!}=\prod_{j=1}^{p-1} (1+\frac{mp}{j})
\equiv 1+mp\sum_{j=1}^{p-1}\frac{1}{j}\pmod {p^2}\equiv 1\pmod{p^2},
\]
since in the last sum $\frac{1}{j}+\frac{1}{p-j}\equiv 0\pmod p$ for $p>2$.
Therefore,
\begin{eqnarray*}
\frac{(c_k p)!}{p^{c_k}}\equiv {c_k}!\, (p-1)!\pmod{p^2}.
\end{eqnarray*}
Taking the product
$\displaystyle\prod_{k=1}^s \frac{(c_k p)!}{p^{c_k}}
\equiv {(p-1)!}^p \prod_{k=1}^s {c_k}! \equiv -\prod_{k=1}^s {c_k}!$,
since ${(p-1)!}^p\equiv -1\pmod p$, which replaced in
\refeq{congr_pq} produces the claim.
 The theorem is proved.
\end{proof}

The following corollary  gives a
 complete answer to the second question
of Hough and Simion.
\begin{cor}
$F(4,n)$ is divisible by $4$, unless
$n$ is of the form $\displaystyle \frac{2^{2t}-1}{3}$, in which case
$F(4,n)\equiv 1\pmod 4$, or of the form
$\displaystyle \frac{2^{2t+1}+2^{2j+1}-1}{3},$
for  $t>j$, in which case $F(4,n)\equiv 2\pmod 4$.
\end{cor}

We include here the base 2 representations of the above numbers
$\leq 60$, namely $1,3$, $5,11$, $13,21$, $43,45,53$
since it suggests a recursive construction of the sequence,
\begin{eqnarray*}
& 1_2, 11_2,\\
& 101_2, 1011_2,1101_2,\\
& 10101_2, 101011_2,101101_2,110101_2.
\end{eqnarray*}
We see that on each row, we start with $1010\ldots 1$,
obtaining the rest of the strings by inserting the bit 1 to the right of
an already existent bit 1, starting with the rightmost one.

What are the possible {\em least} distinct residues of the
$p^2$-Catalan numbers modulo $p^2$?
We give the following table, with the least residues
of $\displaystyle \frac{1}{(p^2-1)n+1}\binom{p^2n}{n}\pmod {p^2}$,
for $p=2,3,5,7$, computed easily by hand, using
Theorem \refth{thm7}. We listed the partitions of $p$, and we computed
the residue modulo $p^2$ of the multinomial coefficient corresponding to
each partition, eliminating duplicates.
For instance, if $p=5$, the partitions of 5 are:
$\{\{ 5\},\{ 4,1\} ,\{ 3,2\}$, $\{ 3,1,1\},\{ 2,2,1\}$,
      $\{ 2,1,1,1\},\{ 1,1,1,1,1\} \}$, so by Theorem \refth{thm7}
the least residues of $F(5^2,n)$ modulo $5^2$ are:
$0$ and $\binom{5}{5}=1$, $\binom{5}{4,1}=5$, $\binom{5}{3,2}=10,$
$\binom{5}{3,1,1}=20$, $\binom{5}{2,2,1}=30\equiv 5\pmod{25}$,
$\binom{5}{2,1,1,1}=60\equiv 10\pmod{25}$, $\binom{5}{1,1,1,1,1}=120\equiv 20\pmod{25}$
\[
\begin{tabular}{|l|l|}
\hline
\strut
 $p$ & least residues modulo $p^2$ of $F(p^2,n)$\strut \\
\hline
\strut
2 & $0,1,2$ \strut \\
\hline
\strut
3 & $0,1,3,6$ \\
\hline
5 &  $0,1,5,10,20$ \\
\hline
\strut
7 & $0,1,7,14,21,35,42$ \\
\hline
\end{tabular}
\]

We provide here the following weak bound.
\begin{cor}
The number of distinct residues of $p^2$-Catalan numbers $\pmod {p^2}$ is
less than or equal to $\pi(p)+1$, where $\pi(p)$ is the number of
partitions of $p$.
\end{cor}
\begin{proof}
Straightforward.
\end{proof}

\begin{rem}
If we denote by $a_s$ the number of distinct residues of
$s^2$-Catalan numbers $\pmod {s^2}$, then $\{a_s\}_s$ is the sequence
{\em A053991} in \cite{OnlineE}.
\end{rem}

\section{Divisibility of $p^q$-Catalan numbers}

Now we attempt to find all natural numbers for which
$\displaystyle p^q$ divides $F(p^q,n)$, $q\geq 3$.
Denote by $j_i$ the least non negative residue of $\al_i\pmod q$.
We prove the result
\begin{thm}
\label{thm10}
When $p$ is an odd prime and $q\geq 3$, then
$p^q$ divides $F(p^q,n)$,
unless $n$ is of the form
$\displaystyle\frac{p^{tq}-1}{p^q-1}$, for some $t\in {\bf N}$,
or of the form
$\displaystyle\frac{p^{qt_1 +j_1}+\cdots +
p^{qt_{m(p-1)+1}+j_{m(p-1)+1}}-1}{p^q-1}$, for some $t_i\in {\bf N}$,
 $1\leq m\leq q-1,\, 0\leq j_i\leq q-1$,  and
$\displaystyle \sum_i p^{j_i}\equiv 1\pmod {p^q-1}$.
\end{thm}

\begin{proof}
By Lemma \refth{pow_lem}, if $p^q\,\not|\, F(p^q,n)$, then
\begin{eqnarray*}
\label{pow_pq}
{v_p}(F(p^q,n))=\frac{\sigma_p((p^q-1)n+1)-1}{p-1}\leq q-1,
\end{eqnarray*}
so
\(
\sigma_p((p^q-1)n+1)\leq (p-1)(q-1)+1.
\)
If $\sigma_p((p^q-1)n+1)=1$, then $(p^q-1)n+1=p^{tq+i}$, for
some $0\leq i\leq q-1$.
Working modulo $p^q-1$ implies $i=0$. Thus,
$\displaystyle n=\frac{p^{tq}-1}{p^q-1}$.
Assume $\sigma_p((p^q-1)n+1)=l, 1<l\leq (p-1)(q-1)+1.$ We obtain
$(p^q-1)n+1=p^{\al_1}+\cdots +p^{\al_l}, \al_1\leq \cdots \leq \al_l$.
Modulo $(p-1)$, this transforms into
\[
(p^q-1)n\equiv -1+l\equiv 0\pmod {p-1},
\]
which will imply $l=m(p-1)+1$, for some $m$. Since $1<l\leq (q-1)(p-1)+1$,
we get $0<m\leq q-1$.
We obtain, for $\al_i=q t_i+j_i$, $0\leq j_i\leq q-1$,
\[
n= \frac{p^{qt_1 +j_1}+\cdots +
p^{qt_{m(p-1)+1}+j_{m(p-1)+1}}-1}{p^q-1},\ m\in {\mathbf N},
\]
with the condition $\displaystyle \sum_i p^{j_i}\equiv 1\pmod {p^q-1}$.
\end{proof}
We use in the next section the following
\begin{cor}
\label{cor}
$p^q\, |\, \frac{1}{(p^q-1)n+1}\binom{p^q n}{n}$
if and only if $p^q\, |\, \binom{p^q n+1}{n}$.
\end{cor}

\section{Squarefree Binomial Coefficients}

In this section we study squarefree binomial coefficients of the form
$\binom{p^q n+1}{n}$, with the help of generalized Catalan numbers.
Thus, in order to study these squarefree binomial coefficients,
it suffices to consider only $n$
 of the form
 $\displaystyle
 \frac{p^{qt_1 +j_1}+\cdots +
p^{qt_{m(p-1)+1}+j_{m(p-1)+1}}-1}{p^q-1},$
$t_i\in {\mathbf N}$, $1\leq m\leq q-1,\, 0\leq j_i\leq q-1$,
such that $\displaystyle \sum_i p^{j_i}\equiv 1\pmod {p^q-1}$.

\def\L{\Lambda}

In \cite{GR}, the authors proved
that if $\displaystyle \binom{n}{k}$ is squarefree,
then $n$ or $n-k$ must be small.
Finding explicit bounds is a much more difficult task.
They showed that $\displaystyle \binom{2n}{n}$ is
squarefree for $n>2^{1617}$,
and used some arguments to simplify the computer's work, in checking
the possible exceptions $n=2^r$.
However, our job is not as hard; we rely on
\cite{GR} and use some estimates on the Chebyshev's
function $\psi(x)=\sum_{d\leq x} \L(d)$, where $\L(d)$ is the
Von Mangoldt's function,  $\L(d)=\log{r}$,
if $d=r^s,\ r$ prime and $\L(d)=0,$ otherwise, to show our results.
Define
 $e(x)=e^{x}$ and
 $\psi(x)=0$, if $x$ is an integer, and $\psi(x)=\{x\}-\frac{1}{2}$,
otherwise, where $\{x\}$ is the fractional part of $x$.

The following lemma proves to be very useful
\begin{lem}
\label{lemI}
If $p^q\leq 99999$, the inequality
\begin{equation}
\label{general_bounds}
\begin{split}
&
\displaystyle 0.9999975\,
\sqrt{p^qn+1}-\displaystyle 1.0000025\,
\sqrt{(p^q-1)n+1}>\\
& 21.683\,
p^{\frac{23q}{48}}\,
n^{\frac{23}{48}}\,
\left(\log{(256 ((p^q-1) n+1))})\right)^{\frac{11}{4}}+
\frac{11}{8}(3\log{n}+2q\log{p} ).
\end{split}
\end{equation}
is true for $n\geq \tau_0$ sufficiently large.
\end{lem}
\begin{proof}
First we prove that $\sqrt{1+x}+\sqrt{1+x-n}\leq 2\sqrt{x},\,
2\leq n\leq x+1$.
By squaring we get
\(
2x+2-n+2\sqrt{(x+1)(x+1-n)}\leq 4x,
\)
which is equivalent to
\(
4(x+1)(x+1-n)\leq (n+2x-2)^2.
\)
The last inequality is equivalent to $n^2-16x+8nx\geq 0$,
which is certainly true if $n\geq 2$.
Now, let $x'=1+x$. We evaluate
\[
\begin{split}
& (1-\alpha)\sqrt{1+x}-(1+\alpha)\sqrt{1+x-n}=
\frac{(1-\alpha)^2 x'-(1+\alpha)^2 (x'-n)}
{(1-\alpha)\sqrt{x'}+(1+\alpha)\sqrt{x'-n}}\\
&=\frac{n(1+\alpha)^2-4 \alpha x'}{(1-\alpha)
\sqrt{x'}+(1+\alpha)\sqrt{x'-n}}\geq \frac{n(1+\alpha)^2-4 \alpha x'}
{(1+\alpha)(\sqrt{x'}+\sqrt{x'-n})}\geq
\frac{n\left(\frac{(1-\alpha)^2}{1+\alpha}-
\frac{4 \alpha}{n(1+\alpha)} x\right)}{2\sqrt{x}}.
\end{split}
\]
Therefore,
\begin{equation}
\label{ineq}
(1-\alpha)\sqrt{1+x}-(1+\alpha)\sqrt{1+x-n}\geq
\left(\frac{(1-\alpha)^2}{1+\alpha}-
\frac{4 \alpha}{1+\alpha}\, \frac{x}{n}
\right) \frac{n}{2\sqrt{x}}.
\end{equation}
Taking $\displaystyle x=p^qn,\, \alpha=\frac{1}{4\cdot 10^5}$,
in \refeq{ineq},
we get
\begin{equation}
\label{rough_approx}
\begin{split}
\displaystyle
& 0.9999975\,
\sqrt{p^qn+1}-\displaystyle 1.0000025\,
\sqrt{(p^q-1)n+1}\\
& \geq
\left(\frac{0.9999975^2}{1.0000025}-
\frac{1}{100000.25}\, p^q
\right)\, \frac{1}{2\sqrt{p^q}}\, n^{\frac{1}{2}}.
\end{split}
\end{equation}
If  $p^q\leq 99999$, then \refeq{rough_approx}
implies our claim that the inequality \refeq{general_bounds} is
true for $n\geq \tau_0$ sufficiently large, since, by \refeq{rough_approx}, the left side is
$O(n^{\frac{1}{2}})$ and the right side is $O(n^{\frac{23}{48}})$.
\end{proof}

Our main result of this section is stated in the next
\begin{thm}
\label{th_gen}
\begin{sloppypar}
Assume $p^q\leq 99999$. Then,  $\binom{p^qn+1}{n}$
is not squarefree for $\displaystyle n\geq\tau_1=\max
\left( \frac{e^{60}-1}{p^q-1}, 5^{10} p^{5q},\tau_0\right)$.
Moreover, the exceptions, for $n<\tau_1$, if they exist, are of the
 form
 $\displaystyle\frac{p^{tq}-1}{p^q-1}$, for any $t\in {\bf N}$,
 or of the form
$\displaystyle\frac{p^{qt_1 +j_1}+\cdots +
p^{qt_{m(p-1)+1}+j_{m(p-1)+1}}-1}{p^q-1}$, for any $t_i\in {\bf N}$,
$1\leq m\leq q-1,\, 0\leq j_i\leq q-1$,  and
$\displaystyle \sum_i p^{j_i}\equiv 1\pmod {p^q-1}$.
\end{sloppypar}
\end{thm}

We proceed to the proof of the theorem.
Let $P=n(p^q n-n+1)(p^q n+1)$.
Corollary 3.2 (p. 82) of \cite{GR} implies
\begin{lem}
Suppose that $\binom{p^qn+1}{n}$ is squarefree. Then,
\begin{equation}
\label{cor3.2}
\begin{split}
& \left|\sum_{d\in I} \psi\left(\frac{p^qn+1}{d}\right)\L(d)\right|+
\left|\sum_{d\in I} \psi\left(\frac{n}{d}\right)\L(d)\right|\\
&+\left|\sum_{d\in I} \psi\left(\frac{(p^q-1)n+1}{d}\right)
\L(d)\right|\geq
\frac{1}{2}
\sum_{{d\in I}, {(d,P)=1}}\L(d),
\end{split}
\end{equation}
where $I$ is the set of integers $d$ in the range
$\sqrt{(p^q-1)n+1}<d\leq \sqrt{p^qn+1}$.
\end{lem}
\begin{sloppypar}
An immediate consequence of Lemma 7.1 of \cite{GR}
(see also \cite{Va}) is
\[
\left|
\sum_{d\in I} \psi\left(\frac{X}{d}\right)\L(d)
\right|
\leq \frac{1}{2R+2}\, \sum_{d\in I}\L(d)+
\left(
\sum_{0<|r|\leq R} |a_r^{\pm}|
\right)
\max_{X\leq x\leq XR}
\left|\sum_{d\in I}
e\left(\frac{x}{d}\right)\L(d) \right|,
\]
where
\(\displaystyle
a_r^{\pm}=\frac{i}{2 \pi (R+1)}\left(
\pi\left(1-\frac{|r|}{R+1}\right)\cot\left(\frac{\pi r}{R+1}\right)+
\frac{|r|}{r}
\right)
\pm \frac{1}{2R+2}\left(1-\frac{|r|}{R+1}\right).
\)
\end{sloppypar}
Taking $R=10$ and using
{\em Mathematica}\footnote{A Trademark of Wolfram Research}
we obtained $\displaystyle \sum_{0<|r|\leq 10}
|a_r^{\pm}|\sim 0.868\leq \frac{86}{99}$,
which implies
\begin{lem}
\[
\left|\sum_{d\in I}\psi\left(\frac{X}{d}\right)\L(d)\right|\leq
\frac{1}{22}\sum_{d\in I} \L(d) +\frac{86}{99}
\max_{X\leq x\leq 10X}\left|\sum_{d\in I}
e\left(\frac{x}{d}\right)\L(d) \right|.
\]
\end{lem}

Using \refeq{cor3.2} and the previous lemma we get
\begin{eqnarray*}
\frac{1}{2} \sum_{
d\in I,\, (d,P)=1} \L(d)
&\leq &
\displaystyle
 \left|\sum_{d\in I}
\psi\left(\frac{p^q n+1}{d} \right)\L(d) \right|+
\left|\sum_{d\in I} \psi\left(\frac{n}{d}\right)\L(d)\right|\\
&+&
\displaystyle
\left|\sum_{d\in I} \psi\left(\frac{(p^q-1)n+1}{d}\right)
\L(d)\right|
\leq
\frac{3}{22}\sum_{d\in I} \L(d)
\\
&+&
\displaystyle
\frac{86}{99}
\max_{p^qn+1\leq x\leq 10(p^qn+1)}\left|\sum_{d\in I}
e\left(\frac{x}{d}\right)\L(d) \right|\\
&+&
\displaystyle
\frac{86}{99}
\max_{(p^q-1)n+1\leq x\leq 10((p^q-1)n+1)}\left|\sum_{d\in I}
e\left(\frac{x}{d}\right)\L(d) \right|
\\
&+&
\frac{86}{99}
\displaystyle
\max_{n\leq x\leq 10n}\left|\sum_{d\in I}
e\left(\frac{x}{d}\right)\L(d) \right|
\\
&\leq &
\frac{3}{22}\sum_{d\in I} \L(d) +\frac{86}{33}
\max_{n\leq x\leq 10(p^qn+1)}\left|\sum_{d\in I}
e\left(\frac{x}{d}\right)\L(d) \right|.
\end{eqnarray*}
Since, $\displaystyle \sum_{d\in I,\, (d,P)>1 }\L(d)
\leq \log{n}+\log{((p^q-1)n+1)}+\log{(p^qn+1)}\leq
3\log{n}+2q\log{p}$, for $n\geq 2$, we obtain
\begin{equation}
\label{part1}
\sum_{d\in I} \L(d)\leq \frac{43}{6} \max_{n\leq x\leq 10(p^qn+1)}\left|
\sum_{d\in I} e\left(\frac{x}{d}\right)\L(d)
\right|
+\frac{11}{8}(3\log{n}+2q\log{p} ).
\end{equation}

Schoenfeld \cite{Sch},  obtained, for $x\geq e^{30}$,
(see also \cite{Per})
\[
\left|\sum_{d\leq x} \L(d)-x\right|< \frac{1}{4 \cdot 10^5}  x.
\]
Since
$\displaystyle\sum_{d\in I} \L(d)=\sum_{d\leq \sqrt{p^qn+1}}\L(d)-
\sum_{d\leq \sqrt{(p^q-1)n+1}}\L(d)$,
we obtain
\begin{equation}
\label{part00}
\begin{split}
\sqrt{p^qn+1}
-&
\displaystyle \frac{1}{4\cdot 10^5} \sqrt{p^qn+1}-\sqrt{(p^q-1)n+1}-
\frac{1}{4\cdot 10^5}  \sqrt{(p^q-1)n+1}\\
-&\frac{11}{8}(3\log{n}+2q\log{p} )=  0.9999975 \sqrt{p^qn+1}-\\
& 1.0000025\sqrt{(p^q-1)n+1}-\frac{11}{8}(3\log{n}+2q\log{p} )\\
& <
\displaystyle \frac{43}{6}
\max_{n\leq x\leq 10(p^qn+1)}\left|
\sum_{d\in I} e\left(\frac{x}{d}\right)\L(d)
\right|,
\end{split}
\end{equation}
for $\displaystyle n\geq \frac{e^{60}-1}{p^q-1}$.

Now, we apply Theorem 9 of \cite{GR}, a consequence of some very important
bounds on exponential sums.
\begin{thm}[Granville-Ramar\'e]
If $k>0$ integer and $y\leq \frac{1}{5} x^{3/5}$, then
\[
\left|
\sum_{y\leq d\leq y'} e\left( \frac{x}{d}\right)\L(d)\right|\leq
\frac{50}{3}\, y\, \left(\frac{x}{y^{\frac{k+3}{2}}}
\right)^{\frac{1}{4(2^k-1)}}
(\log{16y})^{\frac{11}{4}},
\]
for any $y\leq y'\leq 2y$.
\end{thm}

Since $\sqrt{p^qn+1}\leq 2\sqrt{(p^q-1)n+1}$,
using the above theorem of Granville and Ramar\'e we get,
for $n>5^{10} p^{5q}$
 (to have the bound $y\leq \frac{1}{5}x^{3/5}$),
\begin{eqnarray*}
&&\displaystyle\max_{n\leq x\leq 10(p^qn+1)}\left|
\sum_{d\in I} e\left(\frac{x}{d}\right)\L(d) \right| \leq
\max_{n\leq x\leq 10(p^qn+1)} \frac{50}{3}\,  \sqrt{(p^q-1)n+1}\cdot \\
&&
\left(\frac{x}{({(p^q-1)n+1})^{\frac{k+3}{4}}}
\right)^{\frac{1}{4(2^k-1)}}\, \left(\log{(16 \sqrt{(p^q-1)n+1})}
\right)^{\frac{11}{4}}\\
&&=
\frac{50}{3}\, \sqrt{(p^q-1)n+1}\,
\left(\frac{10(p^qn+1)}{({(p^q-1)n+1})^{\frac{k+3}{4}}}
\right)^{\frac{1}{4(2^k-1)}}\cdot \\
&&\left(\log{(16 \sqrt{(p^q-1)n+1})}\right)^{\frac{11}{4}}
\displaystyle \leq \frac{50}{3}\, 11^{\frac{1}{4(2^k-1)}}\,
p^{\frac{q}{4(2^k-1)}}\cdot\\
&& {n}^{\frac{1}{4(2^k-1)}}
\left((p^q-1)n+1\right)^{\frac{1}{2}-\frac{k+3}{4^2(2^k-1)}}\,
2^{-\frac{11}{4}}\,
 \left(\log{(256 (p^q-1)n+1})\right)^{\frac{11}{4}}\leq \\
&&
\frac{50}{3}\, 2^{-\frac{11}{4}}\, 11^{\frac{1}{4(2^k-1)}}\,
p^{q\left(\frac{1}{2}-\frac{k-1}{4^2(2^k-1)}\right) }\,
n^{\frac{1}{2}-\frac{k-1}{4^2(2^k-1)} }\,
 \left(\log{(256 ((p^q-1) n+1))}\right)^{\frac{11}{4}}.
\end{eqnarray*}

We obtain (by taking $k=2$ - that will suffice for our purpose)
\begin{eqnarray*}
\displaystyle
&&\displaystyle\max_{n\leq x\leq 10(p^qn+1)}\left|
\sum_{d\in I} e\left(\frac{x}{d}\right)\L(d) \right| \leq
\frac{50}{3}\, 2^{-\frac{11}{4}}\, 11^{\frac{1}{12}}\,
p^{\frac{23q}{48}}\,
n^{\frac{23}{48} }\,
 \left(\log{(256 ((p^q-1) n+1))}\right)^{\frac{11}{4}}
\end{eqnarray*}

By combining \refeq{part00}
and the previous inequality, we get
\begin{equation}
\label{general_bounds1}
\begin{split}
&
\displaystyle 0.9999975\,
\sqrt{p^qn+1}-\displaystyle 1.0000025\,
\sqrt{(p^q-1)n+1}\leq\\
& 21.683\,
p^{\frac{23q}{48}}\,
n^{\frac{23}{48}}\,
\left(\log{(256 ((p^q-1) n+1))}\right)^{\frac{11}{4}}+
\frac{11}{8}(3\log{n}+2q\log{p} ),
\end{split}
\end{equation}
which is certainly false for $n$ sufficiently large
by Lemma \refth{lemI}. Thus the assumption in Lemma 14 
is false, which implies our Theorem \refth{th_gen}.

\begin{rem}
The inequality \refeq{general_bounds} provides explicit
bounds for $n$ for any choice of $p$ and $q$, with $p^q\leq 99999$.
We can increase the bound for $p^q$, by using a weaker
result of Schoenfeld \cite{Sch}. However, in doing that
we increase the bound on $n$ as well, so we
preferred a better bound on $n$.
\end{rem}

\begin{thm}
$\binom{4n+1}{n}$ is squarefree for
$\displaystyle n\geq 2^{1518}$, and if $n< 2^{1518}$ there are at most
$289,179=\binom{761}{2}\sim 2^{18.2}$ possible exceptions.
\end{thm}
\begin{proof}
If $(p,q)=(2,2)$,
the inequality \refeq{general_bounds}, valid for
$n\geq \max\left(\frac{e^{60}-1}{3},5^{10} 2^{10}  \right)$
changes into
\begin{equation}
\begin{split}
&
\displaystyle 0.9999975\,
\sqrt{4n+1}-\displaystyle 1.0000025\,
\sqrt{3n+1}>\\
& 42.1311\,
n^{\frac{23}{48}}\,
\left(\log{(768 n+1)}\right)^{\frac{11}{4}}+
\frac{33}{8}\log{n}+ 1.65566,
\end{split}
\end{equation}
which is certainly true for $n\geq 2^{1518}$. Theorem \refth{thm7}
and  Corollary \refth{cor} imply that the exceptions (if they exist)
are of the form
$\displaystyle\frac{2^{2t+1}+ 2^{2j+1}-1}{3},\ j\leq t$.
Therefore, since the number of
pairs $(j,t)$, giving different numbers of the above form,
is less than $\binom{761}{2}\sim 2^{18.2}$,
we get the theorem.
\end{proof}

\begin{thm}
$\binom{4n+1}{n}$ is squarefree for
$\displaystyle n\geq 3^{956}$, and if $n<3^{956}$, there are at most
$18088476=\binom{478}{3}\sim 3^{15.3}$ possible exceptions.
\end{thm}
\begin{proof}
If $(p,q)=(3,2)$,
the inequality \refeq{general_bounds}, valid for
$n\geq \max\left(\frac{e^{60}-1}{8},5^{10} 3^{10}  \right)$,
changes into
\begin{equation}
\begin{split}
&
\displaystyle 0.9999975\,
\sqrt{9n+1}-\displaystyle 1.0000025\,
\sqrt{8n+1}>\\
& 26.04\,
n^{\frac{23}{48}}\,
\left(\log{(2048 n+1)}\right)^{\frac{11}{4}}+
\frac{33}{8}\log{n}+ 2.62417,
\end{split}
\end{equation}
which is true for $n\geq 3^{956}$. As in the previous proof,
we get that
the exceptions (if they exist) are of the form
$\displaystyle\frac{2^{2t+1}+ 2^{2j+1}+2^{2i+1}-1}{8},\ i\leq j\leq t$.
Therefore,
since the number of
triples $(i,j,t)$, giving different numbers of the above form,
 is less than $\binom{478}{3}\sim 3^{15.3}$,
we get the result.
\end{proof}

Although, there are not that many possible exceptions,
because of their size, to check each value in acceptable time,
 is beyond any computer's capability at this moment, and
 we were not able to decrease the complexity.
However, we conjecture
\begin{conj}
\label{th2^2}
Except for $1$, $3$ and
$45$, $\binom{4n+1}{n}$ is not squarefree.
\label{th3^2}
Except for $1$, $4$ and $10$,
$\binom{9n+1}{n}$ is not squarefree.
\end{conj}

\vspace{.5cm}

\noindent{\bf Acknowledgements.}
{
The author would like to thank the anonymous referee for his
helpful comments, which improved significantly the presentation
of the paper.
}

\end{document}